\documentclass[10pt,english]{article}
\usepackage[T1]{fontenc}
\usepackage[latin9]{inputenc}
\usepackage[a4paper]{geometry}
\usepackage{color}
\usepackage{babel}
\usepackage{bbding}
\usepackage{url}
\usepackage{amsmath}
\usepackage{amsthm}
\usepackage{amssymb}
\usepackage{setspace}
\usepackage[unicode=true,pdfusetitle,
 bookmarks=true,bookmarksnumbered=false,bookmarksopen=false,
 breaklinks=false,pdfborder={0 0 0},pdfborderstyle={},backref=false,colorlinks=true]
 {hyperref}
\hypersetup{
 pdfencoding=auto, psdextra}

\makeatletter

\providecommand{\LyX}{\texorpdfstring%
  {L\kern-.1667em\lower.25em\hbox{Y}\kern-.125emX\@}
  {LyX}}
\newcommand*{\LyxTextAccent}[3][0ex]{%
  \hmode@bgroup\ooalign{\null#3\crcr\hidewidth
  \raise#1\hbox{#2}\hidewidth}\egroup}
\newcommand{\LyxAccentSize}[1][\sf@size]{%
  \check@mathfonts\fontsize#1\z@\math@fontsfalse\selectfont
}
\ProvideTextCommandDefault{\textcommabelow}[1]{
  \LyxTextAccent[-.31ex]{\LyxAccentSize,}{#1}}

\numberwithin{equation}{section}
\numberwithin{figure}{section}
\theoremstyle{plain}
\newtheorem{thm}{\protect\theoremname}
\theoremstyle{definition}
\newtheorem{defn}[thm]{\protect\definitionname}
\theoremstyle{plain}
\newtheorem{lem}[thm]{\protect\lemmaname}
\theoremstyle{remark}
\newtheorem{rem}[thm]{\protect\remarkname}
\theoremstyle{plain}
\newtheorem{prop}[thm]{\protect\propositionname}

\usepackage{etoolbox}
\usepackage{color}
\usepackage{tikz-cd}
\usepackage[labelformat=empty,belowskip=2pt]{caption}
\usepackage{bbding}
\usepackage{fontawesome}
\usepackage{gauss}

\usepackage{tikz}
\usetikzlibrary{arrows, fit, shapes.geometric}
\tikzset{>=latex}

\usepackage{environ}
\makeatletter
\newsavebox{\measure@tikzpicture}
\NewEnviron{scaletikzpicturetowidth}[1]{%
\def\tikz@width{#1}%
\begin{lrbox}{\measure@tikzpicture}%
\BODY
\end{lrbox}%
\pgfmathparse{#1/\wd\measure@tikzpicture}%
\BODY
}
\makeatother

\tikzset{
seloneEEEEEEE/.style = {scale=.6,fill=magenta!20, inner sep=0pt,outer sep=4pt},
seltwoEEEEEEE/.style = {scale=.6,fill=green!20, inner sep=0pt,outer sep=4pt},
noselEEEEEEE/.style = {scale=.6,inner sep=0pt,outer sep=4pt}
}

\usepackage{graphicx}
\usepackage{mathtools}

\usepackage{breqn}

\newcount \eqinterlinepenalty \eqinterlinepenalty=100

\newcommand{\breakingcomma}{%
  \begingroup\lccode`~=`,
  \lowercase{\endgroup\expandafter\def\expandafter~\expandafter{~\penalty0 }}}

\newcommand{\breakingspace}{%
  \begingroup\lccode`~=32
  \lowercase{\endgroup\expandafter\def\expandafter~\expandafter{~\penalty0 }}}

\usepackage[shadow]{todonotes}

\makeatletter
\patchcmd{\chapter}{\if@openright\cleardoublepage\else\clearpage\fi}{}{}{}
\makeatother

\definecolor{lightgray}{gray}{.80}
\definecolor{lightergray}{gray}{.90}

\newcommand{\elem}[1]{\ifstrequal{#1}{0}{\color{lightgray} 0}{{ #1}}}

\makeatletter
\def\bbordermatrix#1{\begingroup \m@th
  \global\let\perhaps@scriptstyle\scriptstyle
  \@tempdima 4.75\p@
  \setbox\z@\vbox{%
    \def\cr{%
      \crcr
      \noalign{%
        \kern2\p@
        \global\let\cr\endline
        \global\let\perhaps@scriptstyle\relax
      }%
    }%
    \ialign{$\make@scriptstyle{##}$\hfil\kern2\p@\kern\@tempdima
      &\thinspace\hfil$\perhaps@scriptstyle##$\hfil
      &&\quad\ \,\,\hfil$\perhaps@scriptstyle##$\hfil\crcr
      \omit\strut\hfil\crcr
      \noalign{\kern-\baselineskip}%
      #1\crcr\omit\strut\cr}}%
  \setbox\tw@\vbox{\unvcopy\z@\global\setbox\@ne\lastbox}%
  \setbox\tw@\hbox{\unhbox\@ne\unskip\global\setbox\@ne\lastbox}%
  \setbox\tw@\hbox{$\kern\wd\@ne\kern-\@tempdima\left[\kern-\wd\@ne
    \global\setbox\@ne\vbox{\box\@ne\kern2\p@}%
    \vcenter{\kern-\ht\@ne\unvbox\z@\kern-\baselineskip}\,\right]$}%
  \null\;\vbox{\kern\ht\@ne\box\tw@}\endgroup}
\def\make@scriptstyle#1{\vcenter{\hbox{$\scriptstyle#1$}}}
\makeatother

\usepackage{xcolor}

\newcommand{\Mod}{\text{\rm mod-}}
\newcommand{\rep}{\text{\rm rep-}}
\newcommand{\D}{\widetilde{\mathbb{D}}}
\newcommand{\A}{\widetilde{\mathbb{A}}}
\newcommand{\E}{\widetilde{\mathbb{E}}}

\newcommand{\dimz}{\underline\dim}
\newcommand{\End}{\operatorname{End}_{kQ}}
\newcommand{\Hom}{\operatorname{Hom}_{kQ}}
\newcommand{\Homk}{\operatorname{Hom}_k}
\newcommand{\Ext}{\operatorname{Ext}_{kQ}}

\definecolor{darkblue}{HTML}{08009f}
\hypersetup{linkcolor=darkblue}

\DeclareMathSymbol{\shortminus}{\mathbin}{AMSa}{"39}

\makeatletter
 \renewcommand*\l@section{\@dottedtocline{1}{1.5em}{2.3em}}
 \renewcommand*\l@subsection{\@dottedtocline{2}{3.8em}{3.2em}}
 \renewcommand*\l@subsubsection{\@dottedtocline{3}{7.0em}{4.1em}}
\makeatother

\makeatother

\providecommand{\definitionname}{Definition}
\providecommand{\lemmaname}{Lemma}
\providecommand{\propositionname}{Proposition}
\providecommand{\remarkname}{Remark}
\providecommand{\theoremname}{Theorem}

\begin{document}
\title{Proof of the tree module property for exceptional representations
of tame quivers}
\date{\date{}}
\author{Szabolcs Lénárt\thanks{ \Envelope{} \protect\url{lszcs90@gmail.com}, \textbf{Bitdefender
S.R.L.} (400107 Cluj-Napoca, str. Cuza Vod\u{a}, nr. 1, Romania)}, Ábel L\H{o}rinczi\thanks{\Envelope{} \protect\url{lorinczi@math.ubbcluj.ro}, \textbf{Faculty
of Mathematics and Computer Science, Babe\textcommabelow{s}-Bolyai
University} (400084 Cluj-Napoca, str. M. Kog\u{a}lniceanu, nr. 1,
Romania)}, Csaba Szántó\thanks{\Envelope{} \protect\url{szanto.cs@gmail.com}, \textbf{Faculty of
Mathematics and Computer Science, Babe\textcommabelow{s}-Bolyai University}
(400084 Cluj-Napoca, str. M. Kog\u{a}lniceanu, nr. 1, Romania)}, István Szöll\H{o}si\thanks{\Envelope{} \protect\url{szollosi@gmail.com} (corresponding author),
\textbf{Faculty of Mathematics and Computer Science, Babe\textcommabelow{s}-Bolyai
University} (400084 Cluj-Napoca, str. M. Kog\u{a}lniceanu, nr. 1,
Romania), \textbf{Eötvös Loránd University, Faculty of Informatics}
(H-1117 Budapest, Pázmány P. sny 1/C, Hungary)}}
\maketitle
\begin{abstract}
This document serves as an arXiv entry point for the appendix to the
paper \cite{Lenart} (the ancillary file \texttt{e6\_proof.pdf} --
``Proof of the tree module property for exceptional representations
of the quiver $\E_{6}$'') and the appendix to the paper \cite{LLSSd6}
(the ancillary file \texttt{d6\_proof.pdf} -- ``Proof of the tree
module property for exceptional representations of the quiver $\D_{6}$'').
The ancillary files contain the computer generated part of the proofs
of the main results in \cite{Lenart} respectively \cite{LLSSd6},
giving a complete and general list of tree representations corresponding
to exceptional modules over the path algebra of the canonically oriented
Euclidean quiver $\E_{6}$, respectively $\D_{6}$. The proofs (involving
induction and symbolic computation with block matrices) were partially
generated by a purposefully developed computer software, outputting
in a detailed step-by-step fashion as if written ``by hand''.

We also give here a short theoretical introduction and an overview
of the computational method used to prove the formulas given in the
papers \cite{Lenart} and \cite{LLSSd6}.
\end{abstract}

\section{Basic notions of representation theory of algebras\label{sec:Basic-notions}}

Let $Q=(Q_{0},Q_{1},s,t)$ be a\emph{ quiver}, that is, a directed
graph, where $Q_{0}$ is the set of vertices, $Q_{1}$ is the set
of arrows and $s,t:Q_{1}\to Q_{0}$ are functions which attach to
an arrow $\alpha\in Q_{1}$ its source $s(\alpha)\in Q_{0}$ and its
target $t(\alpha)\in Q_{0}$. We often write shortly $Q=(Q_{0},Q_{1})$.
Let $k$ be a field and consider the path algebra $kQ$. The category
$\Mod kQ$ of finite dimensional right modules over $kQ$ can be identified
with the category $\rep kQ$ of the finite dimensional $k$-representations
of the quiver $Q$ (therefore we will use the terms ``module'' and
``representation'' interchangeably).

Recall that a \textbf{$k$}\emph{-representation} $M=(M_{i},M_{\alpha})$
of $Q$ is defined as a set of finite dimensional $k$-spaces $\{M_{i}\,|\,i\in Q_{0}\}$
corresponding to the vertices together with $k$-linear maps $\{M_{\alpha}:M_{s(\alpha)}\to M_{t(\alpha)}\,|\,\alpha\in Q_{1}\}$
corresponding to the arrows. Given two representations $M=(M_{i},M_{\alpha})$
and $N=(N_{i},N_{\alpha})$ of the quiver $Q$, a \emph{morphism of
representations} $f:M\to N$ consists of a family of $k$-linear maps
(corresponding to the vertices) $f_{i}:M_{i}\to N_{i}$, such that
$N_{\alpha}f_{s(\alpha)}=f_{t(\alpha)}M_{\alpha}$ for all $\alpha\in Q_{1}$.
The \emph{dimension vector} of a representation $M=(M_{i},M_{\alpha})$
is 
\[
\dimz M=(d_{i})_{i\in Q_{0}}\in\mathbb{Z}Q_{0}\text{ where }d_{i}=\dim_{k}M_{i},
\]
which is treated as an $n$-dimensional row vector where $n=|Q_{0}|$.
In this case the length of $M$ is $\ell(M)=\sum_{i\in Q_{0}}d_{i}$.

There are five types of so-called \emph{Euclidean (or tame) quivers}:
$\A_{m}$, $\D_{m}$, $\E_{6}$, $\E_{7}$ and $\E_{8}$. 

The \emph{Euler form} of an arbitrary acyclic quiver $Q$ is the bilinear
form defined on $\mathbb{Z}Q_{0}$ as 
\[
\langle x,y\rangle=\sum_{i\in Q_{0}}x_{i}y_{i}-\sum_{\alpha\in Q_{1}}x_{s(\alpha)}y_{t(\alpha)}.
\]
Its quadratic form $q_{Q}$ (called \emph{Tits form}) is independent
from the orientation of $Q$ and in the tame case it is positive semi-definite
with radical $\mathbb{Z}\delta$, where $\delta$ is a minimal positive
imaginary root of the corresponding Kac--Moody root system. A vector
$x\in\mathbb{Z}Q_{0}$ is called \emph{real root} if $q_{Q}(x)=1$,
\emph{imaginary root} if $q_{Q}(x)=0$ and it is \emph{positive} if
$x_{i}\in\mathbb{N}$ for all $i\in Q_{0}$. For two (dimension) vectors
$d,d'\in\mathbb{Z}Q_{0}$ we say that $d\leq d'$ if $d_{i}\leq d'_{i}$
for all $i\in Q_{0}$.

Let $P(i)$ and $I(i)$ be the indecomposable projective respectively
injective module corresponding to the vertex $i$. The \emph{Cartan
matrix} $C_{Q}$ is a matrix with the $j$-th column being equal with
$\underline{\dim}P(j)$. The \emph{Coxeter matrix} is defined as $\Phi_{Q}=-C_{Q}^{t}C_{Q}^{-1}$.
Then $\Phi_{Q}\delta=\delta$ and the Euler form satisfies $\langle a,b\rangle=a\left(C_{Q}^{-1}\right)^{t}b^{t}=-\langle b,\Phi_{Q}a\rangle$,
where $a,b\in\mathbb{Z}Q_{0}$. Moreover, because our algebra is hereditary,
for two modules $M,N\in\Mod kQ$ we get 
\begin{equation}
\langle\dimz M,\dimz N\rangle=\dim_{k}\Hom(M,N)-\dim_{k}\Ext^{1}(M,N).\label{eq:HomExt}
\end{equation}

The \emph{Auslander--Reiten translate}s are defined as 
\[
\tau=D\Ext^{1}(-,kQ)\quad\textrm{and}\quad\tau^{-1}=\Ext^{1}(D(kQ),-)
\]
where $D=\Homk(-,k)$.

An indecomposable module $M$ is \emph{preprojective} if there exists
a positive integer $s$ such that $\tau^{s}(M)=0$, while it is called
\emph{preinjective} if $\tau^{-s}(M)=0$. The indecomposable $M$
is \emph{regular} if it is neither preinjective nor preprojective.

From now on let $Q$ be a tame quiver. For $Q$, the structure of
the category $\Mod kQ$ and its \emph{Auslander--Reiten quiver} is
well-known. Up to isomorphism, the indecomposable preprojective modules
are $\tau^{-s}P(i)$, while the indecomposable preinjectives are $\tau^{s}I(i)$,
where $s\in\mathbb{N}$ and $i\in Q_{0}$. In the sequel we use the
somewhat more convenient notation $P(s,i)$ to denote the preprojective
indecomposable module $\tau^{-s}P(i)$ and $I(s,i)$ to denote the
preinjective indecomposable module $\tau^{s}I(i)$. The following
is true concerning the dimension vectors of preprojective, respectively
preinjective indecomposables: 
\begin{equation}
\dimz P(s,i)=\Phi_{Q}^{-s}\cdot\underline{\dim}P(i)\quad\textrm{and}\quad\dimz I(s,i)=\Phi_{Q}^{s}\cdot\underline{\dim}I(i).\label{eq:dimz}
\end{equation}

The category of regular modules is an abelian, exact subcategory which
decomposes into a direct sum of serial categories with Auslander--Reiten
quiver of the form $\mathbb{Z}\mathbb{A}_{\infty}/r$, called \emph{tube}s
of rank $r$. A tube of rank $1$ is called \emph{homogeneous}, otherwise
it is called \emph{non-homogeneous}.

A very important fact is that $\Mod kQ$ is a \emph{Krull--Schmidt
category}, meaning that every module can be written as a direct sum
of indecomposables in a unique way (up to order and isomorphism).

It is well-known that the dimension vector $x$ of an indecomposable
is either a positive real root (i.e. $q_{Q}(x)=1$) or a positive
imaginary root (i.e. $q_{Q}(x)=0$). It is also known that for every
positive real root $x$ there is a unique (up to isomorphism) indecomposable
$M$ with $\dimz M=x$ (in fact these indecomposables are all the
preprojectives, all the preinjectives and the non-homogeneous regular
indecomposables with dimension different from a multiple of $\delta$). 

An indecomposable module $M$ is called \emph{exceptional}, if it
has no self-extensions (i.e. if $\dim_{k}\Ext^{1}(M,M)=0$). This
means that its dimension is a positive real root (called \textit{exceptional
root}) and $\dim_{k}\operatorname{End}_{kQ}(M)=1$. We know that the
exceptional indecomposable modules are all the preprojectives, all
the preinjectives and the regular non-homogeneous indecomposables
with dimension vector falling below $\delta$ (see \cite{Dlab}). 

For more details concerning the notions presented in this section
we refer to \cite{Aus,Skow1,Skow2,Zhang}. 

\section{Tree representations\label{sec:Tree-representations}}

An indecomposable module $M=(M_{i},M_{\alpha})$ is called a \emph{tree
module} if there is a basis $B$ such that the matrices of the linear
maps $M_{\alpha}$, written in basis $B$ consist only of elements
$0$ and $1$, moreover, the total number of non-zero elements is
$\ell(M)-1$, where $\ell(M)=\sum_{i\in Q_{0}}d_{i}$ with $\underline{\dim}M=(d_{i})_{i\in Q_{0}}$.
Equivalently, $M$ is a tree module if there exists a basis $B$ such
that the associated coefficient quiver is a tree (for details see
\cite{ringel}).

In \cite{ringel} Ringel proves that exceptional modules are tree
modules. The proof is based on a result by Schofield (see \cite{Scho}),
stating that if $M$ is an exceptional module that is not simple,
then there are exceptional modules $X,Y$ with the properties $\Hom(X,Y)=\Hom(Y,X)=\Ext^{1}(Y,X)=0$
and an exact sequence of the following form: $\begin{tikzcd} 0\arrow[r] & vY\arrow[r] & M\arrow[r] & uX\arrow[r] & 0, \end{tikzcd}$
where $u$ and $v$ are positive integers and the notation $uY$ means
$Y\oplus\cdots\oplus Y$ ($u$ times). There are precisely $s(M)-1$
such sequences where $s(M)$ is the number of nonzero components in
$\dimz M$. We call these short exact sequences \emph{Schofield sequence}s
and the pair $(X,Y)$ a \emph{Schofield pair} (associated to $M$).
Note that the original proof of Schofield assumes an algebraically
closed field, but Ringel gives a proof in \cite{Ring2} which works
in arbitrary field $k$. Proposition 6 from \cite{SzSz} states that
if $X$, $Y$, $M$ are exceptional indecomposables such that $u\dimz X+v\dimz Y=\dimz M$,
then we have a Schofield sequence \begin{tikzcd} 0\arrow[r] & vY\arrow[r] & M\arrow[r] & uX\arrow[r] & 0, \end{tikzcd}
if and only if $\left\langle \dimz X,\dimz Y\right\rangle =0$. This
means that Schofield sequences and pairs depend only on the dimensions
of indecomposables, thus their existence condition is field independent.
Also note that although the short exact sequences used in our proofs
are Schofield sequences (as above, with $v=u=1$), we do not use here
the results from \cite{Scho} or \cite{Ring2} to construct them,
but every short exact sequence used throughout the proofs is written
(and verified) using Lemma \ref{lem:SES} (working also over an arbitrary
field $k$).

Although tree representations for some particular quivers are known,
the proof in \cite{ringel} does not give an explicit method for constructing
them in general.

In \cite{Gabriel} Gabriel gave a full list of indecomposable representations
for the Dynkin quivers using $0-1$-matrices. All the given representations
(excepting $4$ of them) were tree representations. Tree representations
in these four cases were given by Crawley-Boevey \cite{Craw2}.

Regarding the Euclidean case, Mróz gave a full list of the indecomposable
tree representations for the quiver of type $\D_{4}$ with four subspace
orientation in \cite{Mroz}. His results were generalized by L\H{o}rinczi
and Szántó, giving a full list of tree representations for the quiver
of type $\D_{6}$ with a particular non-canonical orientation (see
\cite{abel}). We note that these representations were proved for
path algebras over algebraically closed fields only, moreover in the
paper \cite{abel} indecomposability was checked only for some random
representations from the list (so the checking was not complete).
Analogous problems are considered for canonical algebras in \cite{Dow1},
for nilpotent operators in \cite{Dow2} and for poset representations
in \cite{Grz}.

Concerning the $\D_{m}$, and $\E_{8}$ cases, indecomposable representations
for preinjectives and preprojectives were given by Kussin, K\c{e}dzierski
and Meltzer in \cite{Kedz} and \cite{Kussin}, respectively (however,
those representations are not tree representations). Their aim was
not to give explicit tree representations in particular, but to describe
a general method for obtaining indecomposable (not necessarily tree)
representations in tame cases. 

Our first aim was to study tree representations and to develop a computational
method which produces rigorously proved explicit tree formulas (in
a ``ready to consume'' form) and which is also ``scalable'' (can
be performed in a timely manner for all possible families of exceptional
modules). Our second aim was to use the method in producing a complete
and explicit list of tree representations for all families of exceptionals,
which can be easily introduced and used in any computer algebra system,
without bothering about the way they were obtained. Given the nature
of the problem (the number of cases to be considered and the amount
of block-matrix arithmetic to be performed) the best we could come
up with was the idea presented in Subsection 1.3 from \cite{Lenart}
(which, to our knowledge, is new and completely different from the
method(s) used by Mróz, Kussin, K\c{e}dzierski and Meltzer) and to
develop a special proof assistant software performing the matrix-crunching
and producing a rather lengthy, nevertheless completely general, formal
and correct proof of every formula listed in Part II of the ancillary
documents. 

The importance of knowing explicit formulas for tree representations
stems from a number of advantageous properties. In case of tree representations,
the matrices involved are the ``sparsest possible'' (i.e. containing
the minimal number of non-zero elements), thus reducing the storage
and running time complexity in computer implementations. As mentioned
before, the exceptional modules are determined by their dimension
vectors up to isomorphism, so having a formula for each of them gives
a ``nice'' representative of each isomorphism class. In fact, we
could say that tree representations are the ``canonical'' forms
of these modules, analogously to the canonical form of matrix pencils
or canonical forms of matrices (for example the Jordan normal form).
An example of nice consequences of knowing such sparse forms is the
paper \cite{Mroz2} by Mróz, where such matrix forms of modules were
applied to obtain formulas for the multiplicities of the preprojective
and preinjective indecomposables appearing in the decomposition of
an arbitrary $\D_{4}$ module.

It is very important to realize that \emph{the tree representations
given remain valid independently on the underlying field of the representation}.
That is, the $1-0$ matrices listed in the ancillary files withstand
a replacement of the base field $k$ in $\Mod k\Delta(Q)$ such that
if $M\in\Mod k\Delta(Q)$ is an exceptional tree representation, then
$M'\in\Mod k'\Delta(Q)$ is also an exceptional tree representation
where $\dimz M=\dimz M'$, and every matrix $M_{\alpha}$ from the
first representation is formally the same as the corresponding matrix
$M'_{\alpha}$ from the second one.

\section{Proving the field independent tree module property \label{sec:Proving-the-tree}}

In this section we describe the method used to prove the tree module
property for every representation given in the lists in Part II of
the ancillary documents, both from the theoretical and practical perspective.
The method presented here is general (in the sense that it could be
applied to any tame quiver), so as stated before, $Q$ denotes an
arbitrary tame quiver and $k$ an arbitrary field. We just state the
results here, the proofs are to be found in \cite{Lenart}.

We will use the ``field independent'' qualifier in relation to representations
and short exact sequences in the following precise manner:
\begin{defn}
Let $M\in\Mod kQ$ an (exceptional) indecomposable module. We say
that:
\begin{enumerate}
\item[(1)] The module $M$ is \emph{field independent (exceptional) indecomposable}
if in the corresponding representation $M=(M_{i},M_{\alpha})$ all
the elements in the matrices $M_{\alpha}$ are either $0$ or $1$
and for any field $k'$ if we consider a module $M'\in\Mod k'Q$ such
that $\dimz M=\dimz M'$ and every matrix $M'_{\alpha}$ from the
corresponding representation $M'=(M'_{i},M'_{\alpha})$ is formally
the same as $M_{\alpha}$ (for all arrows $\alpha$), then $M'$ is
also (exceptional) indecomposable in $\Mod k'Q$.
\item[(2)] The module $M$ has the \emph{field independent tree property} if
it is a tree module in $\Mod kQ$ and it is also a \emph{field independent}
(exceptional) indecomposable module (i.e. if we consider the corresponding
representation with formally the same matrices over any other field
$k'$, we still get an exceptional indecomposable tree module in $\Mod k'Q$).
\item[(3)] A short exact sequence of the form $\begin{tikzcd} 0\arrow[r] & Y\arrow[r, "f"] & Z\arrow[r, "g"] & X\arrow[r] & 0 \end{tikzcd}$
is \emph{field independent} (with $X,Y,Z\in\Mod kQ$) if all the elements
in the matrices of the representations $X$, $Y$ and $Z$ are either
$0$ or $1$, all the elements in the matrices $f_{i}$ and $g_{i}$
of the embedding $f=(f_{i})_{i\in Q_{0}}$ respectively the projection
$g=(g_{i})_{i\in Q_{0}}$ are either $0$ or $1$ or $-1$ and in
any field $k'$ the sequence $\begin{tikzcd} 0\arrow[r] & Y'\arrow[r, "f'"] & Z'\arrow[r, "g'"] & X'\arrow[r] & 0 \end{tikzcd}$
is also exact, where $X',Y',Z'\in\Mod k'Q$, $f':Y'\to Z'$, $g':Z'\to X'$
correspond in order to $X$, $Y$, $Z$, $f:Y\to Z$, $g:Z\to X$
with the respective dimension vectors unchanged and with all matrices
(both from the representations and from the morphisms) being formally
the same when considering them over $k'$ instead of $k$.
\end{enumerate}
\end{defn}

The following proposition and lemmas constitute the theoretical elements
of the technique used to prove the formulas in a field independent
way:
\begin{lem}
\label{lem:DimEnd1}For a module $M\in\Mod kQ$ we have $M$ is exceptional
indecomposable if and only if $\dim_{k}\End(M)=1$ and $\dimz M\neq\delta$.
\end{lem}

\begin{rem}
\label{rem:DimEnd1}We know that exactly these are the exceptional
modules in the tame case: the preprojective indecomposables, the preinjective
indecomposables and the regular non-homogeneous indecomposable modules
with dimension vector falling below $\delta$.
\end{rem}

\begin{prop}
\label{prop:Two-SES}Let $X,Y,X',Y'\in\Mod kQ$ be indecomposable
modules. If $M\in\Mod kQ$ such that 
\begin{enumerate}
\item[(a)] there is an exceptional $Z\in\Mod kQ$ such that $(X,Y)$ and $(X',Y')$
are Schofield pairs associated to $Z$,
\item[(b)] there exist two short exact sequences \[\begin{tikzcd} 0\arrow[r] & Y\arrow[r] & M\arrow[r] & X\arrow[r] & 0 \end{tikzcd}\]
and \[\begin{tikzcd} 0\arrow[r] & Y'\arrow[r] & M\arrow[r] & X'\arrow[r] & 0, \end{tikzcd}\] 
\item[(c)] $X\ncong X'$ or $Y\ncong Y'$, 
\item[(d)] $\dim_{k}\Ext^{1}(X,Y)=\dim_{k}\Ext^{1}(X',Y')=1$ 
\end{enumerate}
then $M$ is exceptional indecomposable.
\end{prop}

\begin{lem}
\label{lem:SES}Let $X,Y,Z\in\Mod kQ$ and $f=(f_{i})_{i\in Q_{0}}$,
$g=(g_{i})_{i\in Q_{0}}$ families of $k$-linear maps $f_{i}:Y_{i}\to Z_{i}$,
\textup{$g_{i}:Z_{i}\to X_{i}$.} Then there is a short exact sequence
$$\begin{tikzcd} 0\arrow[r] & Y\arrow[r, "f"] & Z\arrow[r, "g"] & X\arrow[r] & 0 \end{tikzcd}$$
if and only if the following conditions hold (we identify the maps
$f_{i}$ and $g_{i}$ with their matrices in the canonical basis):
\begin{enumerate}
\item[(a)] the matrices $f_{i}$ (respectively $g_{i}$) have maximal column
(respectively row) ranks,
\item[(b)] $f_{t(\alpha)}Y_{\alpha}=Z_{\alpha}f_{s(\alpha)}$ and $g_{t(\alpha)}Z_{\alpha}=X_{\alpha}g_{s(\alpha)}$,
for all $\alpha\in Q_{1}$,
\item[(c)] $g_{i}f_{i}=0$, for all $i\in Q_{0}$,
\item[(d)] $\dimz Z=\dimz X+\dimz Y$.
\end{enumerate}
\end{lem}

The embedding $f:Y\to Z$ can be given via a family of maximal (column)
rank matrices $f_{i}$ ($i\in Q_{0}$) satisfying $f_{t(\alpha)}Y_{\alpha}=Z_{\alpha}f_{s(\alpha)}$
for all $\alpha\in Q_{1}$, while the projection $g:Z\to X$ can be
given via a family of maximal (row) rank matrices $g_{i}$ ($i\in Q_{0}$)
satisfying $g_{t(\alpha)}Z_{\alpha}=X_{\alpha}g_{s(\alpha)}$ for
all $\alpha\in Q_{1}$. 
\begin{lem}
\label{lem:Ext1}If $X,Y\in\Mod kQ$ are indecomposable modules such
that $X$ is regular and $Y$ is preprojective, or $X$ is preinjective
and $Y$ is regular or both of them are preprojectives (or preinjectives)
and there is a path in the Auslander--Reiten quiver from the vertex
corresponding to $Y$ to the vertex corresponding to $X$, then $\dim_{k}\Ext^{1}(X,Y)=-\langle\dimz X,\dimz Y\rangle$.
\end{lem}

We are now ready to describe the process of proving the formulas from
Part II of the ancillary document.

\subsection*{The process of proving the field independent tree property}

Suppose we have formulas defining families of matrices $(M_{\alpha}^{(n)})_{\alpha\in Q_{1}}$
depending on some $n\in\mathbb{N}$. The elements of the matrices
$M_{\alpha}^{(n)}$ are either $0$ or $1$, so they can be considered
over an arbitrary field $k$. We want to prove that the representation
of the quiver $Q$ given as $M=M^{(n)}=(M_{i}^{(n)},M_{\alpha}^{(n)})$
has the field independent tree property (where the dimension of each
$k$-space $M_{i}^{(n)}$ is in accordance with the column and row
sizes of the matrices $M_{\alpha}^{(n)}$, thus the formulas also
determine $\dimz M$). Suppose that $\dimz M$ is such that it coincides
with the dimension vector of an exceptional indecomposable (see Lemma
\ref{lem:DimEnd1} and Remark \ref{rem:DimEnd1}). Suppose also that
the number of elements equal to $1$ in the matrices $M_{\alpha}^{(n)}$
is exactly $\ell(M)-1$. So, in order to prove the field independent
tree module property, we need only to show that $M$ is field independent
indecomposable. We may use one of the following lines of reasoning:
\begin{enumerate}
\item[(1)] \emph{Prove that $\dim_{k}\End(M)=1$ in any field $k$ and use Lemma
\ref{lem:DimEnd1}}. This may be done by writing the matrix $A$ of
the homogeneous system of linear equations defining $\End(M)$ and
showing that the corank of $A$ is one (i.e. the solution space is
one dimensional). In order to compute the rank of $A$, it must be
echelonized (brought to row echelon form) using elementary operations
on rows and/or columns in a ``field independent way''. This means
that every single elementary operation used in the process of echelonizing
$A$ must be such that the elements in the resulting matrix are either
$0$, $1$ or $-1$ and the result is exactly the same if performed
in any field $k$. For example if in the case of the matrix $\begin{bmatrix}1 & -1\\
1 & 1
\end{bmatrix}$ we perform the elementary row operation $r_{2}\leftarrow r_{2}-r_{1}$,
then we get $\begin{bmatrix}1 & -1\\
1 & 1
\end{bmatrix}\xrightarrow{r_{2}\leftarrow r_{2}-r_{1}}\begin{bmatrix}1 & -1\\
0 & 2
\end{bmatrix}$ if performed in $\mathbb{R}$, or $\begin{bmatrix}1 & -1\\
1 & 1
\end{bmatrix}\xrightarrow{r_{2}\leftarrow r_{2}-r_{1}}\begin{bmatrix}1 & -1\\
0 & 0
\end{bmatrix}$ if performed in $\mathbb{Z}_{2}$. Hence it has different ranks if
considered over different fields. A crucial element of this proof
is to ensure something like this never happens, but the result of
every single elementary operation performed is formally the same matrix,
independently of the field it is considered in.
\item[(2)] \emph{Perform an induction on $n$, making use of Proposition \ref{prop:Two-SES}}.
First prove the formula for the starting values of $n$ using method
(1) above (typically for $n=0$, but the structure of the block matrices
depending on $n$ might require to make additional proofs for small
values of $n$). Then suppose the formula gives field independent
exceptional indecomposables $M^{(n')}=(M_{i}^{(n')},M_{\alpha}^{(n')})$
for all $n'<n$. Find two pairs of modules $(X,Y)$ and $(X',Y')$
conforming to all requirements of Proposition \ref{prop:Two-SES},
such that any of these four representations is obtained either using
formula $M^{(n')}$ for some $n'<n$ (or some permuted version of
it) or some other formulas proved already to give field independent
exceptional indecomposables. If the quiver $Q$ presents some symmetries,
then a permuted version of the formula $\tilde{M}^{(n')}=(\tilde{M}_{i}^{(n')},\tilde{M}_{i\to j}^{(n')})$
may also be used in the induction step, where $(\tilde{M}_{i}^{(n')})_{i\in Q_{0}}=(M_{\sigma(i)}^{(n')})_{i\in Q_{0}}$
and $(\tilde{M}_{i\to j}^{(n')})_{(i\to j)\in Q_{1}}=(M_{\sigma(i)\to\sigma(j)}^{(n')})_{(i\to j)\in Q_{1}}$
for some permutation $\sigma$. One has to construct here the two
field independent short exact sequences of the form $0\to Y\to M^{(n)}\to X\to0$
and $0\to Y'\to M^{(n)}\to X'\to0$ in order to show their existence.
Once the matrices of the morphisms are constructed, Lemma \ref{lem:SES}
can be used to prove that indeed these form short exact sequences
in any field $k$. We emphasize that conditions (a), (b) and (c) from
Lemma \ref{lem:SES} must be verified in a ``field independent way'':
the rank of the matrices must be checked using field independent echelonization
as explained before, and the result of the matrix arithmetic operations
used in (b) and (c) must be formally the same, independently of the
underlying field.
\item[(3)] \emph{Perform a direct proof, making use of Proposition \ref{prop:Two-SES}}.
Use two pairs of modules $(X,Y)$ and $(X',Y')$ conforming to all
requirements of Proposition \ref{prop:Two-SES}, such that any of
these four representations are obtained by some formulas showed already
to give field independent exceptional indecomposables, and prove the
existence of the two field independent short exact sequences $0\to Y\to M^{(n)}\to X\to0$
and $0\to Y'\to M^{(n)}\to X'\to0$ by constructing them using Lemma
\ref{lem:SES} in the ``field independent way''.
\end{enumerate}
\begin{rem}
\label{rem:Ext1}Note that in methods (2) and (3) the condition $\dim_{k}\Ext^{1}(X,Y)=\dim_{k}\Ext^{1}(X',Y')=1$
required by (d) from Proposition \ref{prop:Two-SES} may be checked
by simply computing $-\langle\dimz X,\dimz Y\rangle$ and $-\langle\dimz X',\dimz Y'\rangle$,
if both pairs are such that Lemma \ref{lem:Ext1} may be applied in
their case.
\end{rem}

The proof process described is extremely cumbersome, time-consuming
and error-prone if performed by a human, therefore we have implemented
a proof assistant software to help us in carrying it out. The proof
assistant can perform any of the steps (1), (2) or (3) based on some
input given in a \LaTeX{} file. The input data consists of the formulas
$(M_{\alpha}^{(n)})_{\alpha\in Q_{1}}$ defining the representations
and the choice for the short exact sequences required in (2) and (3),
together with the families of matrices defining the morphisms. All
this data must be given in a \LaTeX{} document with a well-defined
structure, in order for the proof assistant to be able to parse it
and extract the relevant information. The matrices are given either
as ``usual matrices'' (of fixed size, with elements equal to either
$1$, $-1$ or $0$), or symbolic block-matrices of variable size,
depending on the parameter $n\in\mathbb{N}$. Every block-matrix is
built using the following three types of blocks: zero block of size
$n_{1}\times n_{2}$, the identity block $I_{n}$ and a block denoted
by $E_{n}$ having ones on the secondary diagonal and zeros everywhere
else (note that $E_{n}^{2}=I_{n}$ in every field). We have used the
document processor \LyX{} to edit the input document and export it
to \LaTeX{} (in this way ensuring a syntactically correct \LaTeX{} file).

These are the steps performed by the software:
\begin{itemize}
\item It reads and stores the data $M^{(n)}=(M_{i}^{(n)},M_{\alpha}^{(n)})$
defining the representations $M^{(n)}$.
\item Computes the total number of elements equal to $1$ in the matrices
$M_{\alpha}^{(n)}$ and compares it against $\ell(M^{(n)})$ to ensure
their number is exactly $\ell(M^{(n)})-1$.
\item If instructed to perform along method (1), it computes the matrix
$A$ of the homogeneous system of linear equations defining $\End(M^{(n)})$
and shows that it can be brought to echelon form by performing exactly
the same elementary operations resulting in exactly the same matrix
(formally) if considered in any field. In this way it ensures that
the corank of $A$ is one independently of the field. Note that it
can perform in this mode only with formulas where $n$ has any given
concrete value. 
\item If instructed (and given sufficient data) it performs all checks required
by methods (2) or (3) based on Proposition \ref{prop:Two-SES}. First
it checks in the list provided in \cite{SzSz} to see that both pairs
$(X,Y)$ and $(X',Y')$ are Schofield pairs associated to $Z\in\Mod kQ$
exceptional indecomposable such that $\dimz Z=\dimz M^{(n)}$, then
verifies conditions (c) and (d) from Proposition \ref{prop:Two-SES}.
It is ensured that the requirements of Lemma \ref{lem:Ext1} are met
and condition (d) is validated as mentioned in Remark \ref{rem:Ext1}.
Finally, it ensures the existence of two short exact sequences of
the form $\begin{tikzcd} 0\arrow[r] & Y\arrow[r, "f"] & M^{(n)}\arrow[r, "g"] & X\arrow[r] & 0 \end{tikzcd}$
and $\begin{tikzcd} 0\arrow[r] & Y'\arrow[r, "f'"] & M^{(n)}\arrow[r, "g'"] & X'\arrow[r] & 0 \end{tikzcd}$
by reading the matrices of the morphisms $f$, $f'$, $g$ and $g'$
and showing that every elementary operation and block-matrix arithmetic
may be performed in a field independent way in order to fulfill every
requirement of Lemma \ref{lem:SES}.
\end{itemize}
Every single operation performed by the proof assistant software is
written to this output \LaTeX{} document. Everything (including the
elementary operations and the details of computing the block matrix
sums and products) is output a detailed step-by-step fashion as if
written ``by hand''. In this way one does not have to believe in
the correctness of the implementation, because the complete proof
is ``on paper'' and every single step may be crosschecked and verified
by a human mathematician. 

\section{About this document\label{sec:About-this-document}}

The purpose of this document is to give an overview of the computational
method used to prove the formulas given in the papers \cite{Lenart}
and \cite{LLSSd6} and also to serve as an entry point on arXiv to
the quite lengthy proofs given as separate files. The documents containing
the complete proofs have considerable sizes, so they are given as
attached ancillary documents:
\begin{itemize}
\item the file named \texttt{e6\_proof.pdf} has the title ``Proof of the
tree module property for exceptional representations of the quiver
$\E_{6}$'' and is the appendix to the paper \cite{Lenart};
\item the file named \texttt{d6\_proof.pdf} has the title ``Proof of the
tree module property for exceptional representations of the quiver
$\D_{6}$'' and is the appendix to the paper \cite{LLSSd6}.
\end{itemize}
The ancillary files contain the output generated by the proof assistant
software. Being relatively self-contained materials, the introductory
text from the current document is also present in the appendices.

\end{document}